\documentclass[11pt,a4paper]{article}

\makeatletter \@addtoreset{equation}{section} \makeatother

\usepackage[centertags]{amsmath}
\usepackage{amssymb}
\usepackage{amsfonts}
\usepackage{epsfig}
\usepackage{amssymb}
\usepackage{newlfont}
\usepackage{graphicx}
\usepackage{epstopdf}
\usepackage{cite}
\usepackage{caption}
\usepackage{subcaption}
\usepackage{titlesec}
\usepackage{url}
\usepackage{algpseudocode}
\usepackage{algorithm}
\usepackage{arydshln,leftidx,mathtools}
\usepackage{tikz}

\usetikzlibrary{arrows,matrix,positioning}

\newtheorem{thm}{Theorem}[section]
\newcommand{\be}{\begin{equation}}
\newcommand{\ee}{\end{equation}}
\newcommand{\bes}{\begin{equation*}}
\newcommand{\ees}{\end{equation*}}
\newcommand{\beqa}{\begin{eqnarray*}}
\newcommand{\eeqa}{\end{eqnarray*}}

\renewcommand{\vec}{\mathbf}
\newcommand{\rank}{\mathrm{rank}}
\newcommand{\R}{\mathbb R}

\newcommand{\tha}{\theta}
\newcommand{\shi}{\tau}
\newcommand{\om}{\omega}

\newcommand{\qed}{\hfill \ensuremath{\Box}}

\begin{document}
\title{Linear least squares problems involving fixed knots polynomial splines and their singularity study}

\author{Zahra Roshan Zamir and Nadezda Sukhorukova,\\
Swinburne University of Technology\\
zroshanzamir@swin.edu.au and nsukhorukova@swin.edu.au\\
Swinburne University of Technology, FSET, Department of Mathematics\\
Mail Services Unit, H55\\
PO Box 218\\
HAWTHORN VIC 3122 \\
}



  \maketitle 
  \begin{abstract} In this paper, we study a class of approximation problems,  appearing in data approximation and signal processing. The approximations are constructed as combinations of polynomial splines (piecewise polynomials), whose parameters are subject to optimisation, and so called prototype functions, whose choice is based on the application, rather than optimisation. The corresponding optimisation problems can be formulated as Linear Least Squares Problems (LLSPs). If the system matrix is non-singular, then the corresponding problem can be solved using the normal equations method, while for singular cases slower (but more robust) methods have to be used.  In this paper we develop sufficient conditions for non-singularity. These conditions can be verified much faster than the direct singularity verification of  the system matrices. Therefore, the algorithm efficiency can be improved by choosing a better suited method for solving the corresponding LLSPs.
  \end{abstract}
  
   
\section{Introduction}\label{sec:introduction}
In this paper, we consider two optimisation problems  (Model~1 and Model~2), frequently appearing in approximation and signal processing. The signal approximations are constructed as  products of two functions.  In the  corresponding function products, one of the functions is a polynomial spline (piecewise polynomial) and another function is a prototype (also called basis) function, defined by the application. Common examples of prototype functions are sine and cosine functions. In Model~1 the wave is oscillating around ``zero level'', while Model~2   admits a vertical shift. In the case when the vertical shift is required (Model~2), we construct it as another polynomial spline. 

The choice of polynomial splines is due to the fact that these functions combine the simplicity of polynomials and additional approximation flexibility, which can be achieved by switching from one polynomial to another. Therefore, on the one hand, the corresponding optimisation problems can be solved efficiently and, on the other hand, the approximation inaccuracy (evaluated as the sum of the corresponding deviation squares) is low. 

It is shown  in this paper that Model~1 and Model~2 can be formulated as Linear Least Squares Problems (LLSPs). There are a variety of methods to solve such problems~\cite{GOLUB, LH, AB}. If it is known that the corresponding system matrices are non-singular, then the most popular approach for solving the corresponding LLSPs is based on the normal equations method, since this method is very efficient (fast and accurate) when working with non-singular matrices. If the corresponding matrices are singular, one needs to apply more robust methods,  for example, QR decomposition or Singular Value Decomposition (SVD). These methods are much more computationally expensive.

In this paper we provide sufficient conditions for non-singularity (both models). The obtained conditions are easier to check than the direct singularity verification of  the corresponding matrices. Therefore, the algorithm efficiency can be enhanced by choosing a better suited approach for solving the corresponding LLSPs.  

In this paper, we use truncated power basis function to define splines. Another way to construct basis function is through B-splines. B-splines have several computational advantages when running numerical experiments (have smallest possible support, see~\cite{NUG}). However, truncated power functions are very common when theoretical properties of the models are the research targets (see, for example~\cite{NUG}).

The paper is organised as follows. In section~\ref{sec:optimisation_models} we formulate the optimisation problems (both models). Then, in section~\ref{sec:singularity_study} we develop sufficient conditions for non-singularity. In section~\ref{sec:applications} we provide a detailed example of how our conditions can be applied. Finally, in section~\ref{sec:conclusions} we summarise the obtained results and indicate further research directions.

\section{Approximation models}\label{sec:optimisation_models}
In this section, we formally introduce polynomial splines and explain why these functions are used in our models. Then we formulate the models as mathematical programming problems and demonstrate that this problems are LLSPs.
\subsection{Polynomial splines}\label{subsec:polynomial_splines}
Polynomial splines (piecewise polynomials) combine the simplicity of polynomials and the flexibility that enables them to change abruptly at the points of switching from one polynomial to another (spline knots). These special properties are essential for good quality approximation, since, on the one hand, the corresponding optimisation problems are relatively inexpensive to solve and, on the other hand, the obtained approximations are accurate enough for reflecting the essential characteristics of the original signal (raw data). Therefore, polynomial splines are very commonly used in approximation~\cite{NUG, Meinardus89, Schumaker68}.

There are many ways to construct polynomial splines. One possibility is to do it through truncated power functions~\cite{NUG}:
\begin{equation}\label{eqn4}
S_{m} (\vec x\,,\tha\,, t)= x_{0}+ \sum_{j=1}^m x_{1j
}t^j+\sum_{l=2}^n\sum_{j=1}^m x_{ij}((t-\tha_{l-1})_{+})^j\,,
\end{equation}
where $m$ is the degree,
$\tha=(\theta_{1},\dots,\theta_{n-1})$ are the knots, $x_{0},
x_{11},\dots,x_{nm}$ are the spline parameters and
\begin{equation*}
(t-\tha_{l-1})_{+}=\max\{0,(t-\tha_{l-1})\}=\left\{
\begin{array}{rl}
t-\tha_{l-1}\,, & \text{if}\quad t>\tha_{l-1}\,,\\
0\,, &\text{if}\quad  t\le\tha_{l-1}\,,
\end{array}
\right.
\end{equation*}
is the truncated power function. Note that this presentation implies that the polynomial splines are continuous. In general, polynomial splines can be discontinuous at their knots~\cite{NUG}. In this paper, however, we only consider continuous polynomial splines.

The spline knots can be free or fixed. If the knots are free then they are considered as additional variables in the corresponding optimisation problem  and thus $\vec x=
(x_{0}, x_{11},\dots,x_{nm},\tha)$ (see next section for details). In this case, the corresponding optimisation problem becomes more complex~\cite{NUG, su10, archive, nurnberger892, FreeKnotsOpenProblem96}. In particular, it becomes non-convex. Generally, it is much easier to solve a higher dimension
fixed knots problem (with a considerably larger number of intervals) than a free knots one.

\subsection{Optimisation problems}\label{subsec:optimisation_models}

In this section we formulate our models as mathematical programming problems. We consider two types of models. Model~1 corresponds to the case when the wave is oscillating around ``zero level'', while Model~2 enables a vertical shift in a form of a polynomial spline. Similar problems have been considered in~\cite{ANZ_us, Four_Models}, where it appeared that Model~2, in general, is more computationally expensive, but its approximation accuracy is significantly higher.  

Consider a signal segment $y=(y_1,\dots,y_N)\in\R^N$, where $y_i,~i=1,\dots,N$ are evaluated at the time moments $t_i,~i=1,\dots,N.$  This signal is to be approximated by a function $f(t)$ from one of the following two models.
\begin{equation}\label{eq:approximation1}
\text{Model~1}:~f(t)=A(M,t)g(t)
\end{equation} 
and
\begin{equation}\label{eq:approximation2}
\text{Model~2}:~f(t)=A_1(M,t)g(t)+A_2(L,t),
\end{equation}
where $A(M,t),~A_1(M,t)$ and $A_2(L,t)$ are  polynomial splines with fixed knots, $A, A_1$ and $A_2$ are the corresponding spline parameters, and $g(t)$ is a prototype  function (e.g., sine, cosine). In most cases, the choice of the prototype functions is application driven (specified by application experts, e.g., engineer).

These two approximation problems can be formulated as mathematical programming problems.
\be\label{eq:model1}
\text{Model~1}:~\min_{\vec x} \sum_{i=1}^N (y_{i}-A(\vec x ,t_i)g(t_i))^2\,
\ee
and
\be\label{eq:model2}
\text{Model~2}:~\min_{\vec{x_1},\vec{x_2}} \sum_{i=1}^N (y_{i}-A_1(\vec x_1,t_i)g(t_i)-A_2(\vec x_2,t_i))^2\,
\ee

In the following two sections we study these models in depth.

\subsubsection{Linear Least Squares Optimisation Model 1} 
\label{subsub:LLSOM1}
Assume that the spline degree is $m$, the number of subintervals is $n$ and the corresponding knots are 
$$\theta_0=t_1\leq\theta_1\leq\theta_2\leq\dots\leq\theta_{n-1}\leq\theta_n=t_N.$$
Model~1 is an LLSP since it can be rewritten as follows:
\be\label{eqn7}
\min_{\vec x} \|B\vec x-y\|_{2}^2\,,
\ee
where $y=(y_1,\dots,y_N)^T\in\R^N,$ $y_{i},~i=1,\dots,N$ are the recorded signals at $t_{i},~i=1,\dots,N$ and 
\be\label{eq:matB}
B=\left[
\begin{array}{ccccccccccc}
\alpha_{1} & \alpha_{1}t_{1} & \dots & \alpha_{1}t_{1}^m &  \alpha_{1}\beta_{11} &\dots & \alpha_{1}\beta_{11}^m &\dots & \alpha_{1}\beta_{1n-1} &\dots & \alpha_{1}\beta_{1n-1}^m\\
\alpha_{2} & \alpha_{2}t_{2} & \dots &\alpha_{2}t_{2}^m & \alpha_{2}\beta_{21}& \dots &\alpha_{2}\beta_{21}^m & \dots & \alpha_{2}\beta_{2m}& \dots & \alpha_{2}\beta_{2n-1}^m \\
\vdots&\vdots&\ddots&\vdots&\vdots&\ddots&\vdots&\ddots&\vdots&\ddots&\vdots\\
\alpha_{N} & \alpha_{N}t_{N} & \dots &\alpha_{N}t_{N}^m & \alpha_{N}\beta_{N1}& \dots & \alpha_{N}\beta_{N1}^m & \dots & \alpha_{N}\beta_{Nm}& \dots & \alpha_{N}\beta_{Nn-1}^m
\end{array}
\right ]
\,,\ee
where  $$\beta_{ij}=\max\{0,t_i-\theta_j\},~j=1,\dots,n-1,~i=1,\dots,N\,,$$
and $\alpha_i=g(x_i),~i=1,\dots,N\,.$ 
There exist various methods for solving an LLSP. Most of them are based on the normal equations method, QR decomposition and SVD (see \cite{SL, ABLSP, TB} for details). If $B\in \R^{N\times (mn+1)}$ is a full-rank matrix then, the corresponding LLSP can be solved through the system of normal equations:
 \begin{equation}\label{eq:normal_eq}
 (B^TB)\vec x=B^Ty,
 \end{equation}
where $y=(y_1,\dots,y_N)^T\in \R^N$ is a signal segment recorded at $N$ distinct consecutive time moments. This method is much faster than QR decomposition or SVD but not so accurate.when  matrix $B$ is singular. Therefore, it is essential to develop a singularity testing procedure to choose a suitable method for solving LLSPs.

\subsubsection{Linear Least Squares Optimisation Model 2}\label{subsub:LLSOM2}
In this model, we assume that the wave (signal) is shifted vertically by a spline function. Similar to Model~1, the spline degree is $m$, the number of subintervals is $n$ and the corresponding knots are 
$$\theta_0=t_1\leq\theta_1\leq\theta_2\leq\dots\leq\theta_{n-1}\leq\theta_n=t_N.$$
The corresponding optimisation problem is an LLSP,  formulated as follows:
\be\label{eqn8}
\min_{\vec x} \|B\vec x-y\|_{2}^2\,,
\ee
where $\vec x=[\vec{x_1},\vec{x_2}],$ $y\in \R^N$ is the original signal (see (\ref{eq:model2}) for details) and  $B\in \R^{N\times (2mn+2)}$. Matrix $B$ can be constructed as  
$$B^{N \times (2mn+2)}=[B_{1}^{N \times (mn+1)} \quad  B_{2}^{N \times (mn+1)}],$$ where $B_{1}^{N\times (mn+1)}$ and $B_{2}^{N \times (mn+1)}$ are as follows:
\be\label{matv1}
B_{1}=\left[
\begin{array}{ccccccccccc}
\alpha_{1} & \alpha_{1}t_{1} & \dots & \alpha_{1}t_{1}^m &  \alpha_{1}\beta_{11} &\dots & \alpha_{1}\beta_{11}^m &\dots & \alpha_{1}\beta_{1n-1} &\dots & \alpha_{1}\beta_{1n-1}^m\\
\alpha_{2} & \alpha_{2}t_{2} & \dots &\alpha_{2}t_{2}^m & \alpha_{2}\beta_{21}& \dots &\alpha_{2}\beta_{21}^m & \dots & \alpha_{2}\beta_{2n-1}& \dots & \alpha_{2}\beta_{2n-1}^m \\
\vdots&\vdots&\ddots&\vdots&\vdots&\ddots&\vdots&\ddots&\vdots&\ddots&\vdots\\
\alpha_{N} & \alpha_{N}t_{N} & \dots &\alpha_{N}t_{N}^m & \alpha_{N}\beta_{N1}& \dots & \alpha_{N}\beta_{N1}^m & \dots & \alpha_{N}\beta_{Nn-1}& \dots & \alpha_{N}\beta_{Nn-1}^m
\end{array}
\right ]
\,,\ee
and
\be\label{matv2}
B_{2}=\left[
\begin{array}{ccccccccccc}
1 & t_{1} & \dots & t_{1}^m & \beta_{11} &\dots &  \beta_{11}^m &\dots & \beta_{1n-1} &\dots & \beta_{1n-1}^m\\
1 & t_{2} & \dots & t_{2}^m & \beta_{21}& \dots &\beta_{21}^m & \dots & \beta_{2n-1}& \dots &\beta_{2n-1}^m \\
\vdots&\vdots&\ddots&\vdots&\vdots &\ddots & \vdots & \ddots & \vdots &\ddots & \vdots\\
1& t_{N} & \dots & t_{N}^m & \beta_{N1}& \dots & \beta_{N1}^m & \dots & \beta_{Nn-1}& \dots & \beta_{Nn-1}^m
\end{array}
\right ]
,\ee
where  
$$\beta_{ij}=\max\{0,t_i-\theta_j\},~j=1,\dots,n-1,~i=1,\dots,N$$
and
$\alpha_i=g(t_i),~i=1,\dots,N.$ Note that $B_1$ is similar to $B$ described in (\ref{eq:matB}). In addition, the rows of $B_1$ are obtained from the rows of $B_2$ through multiplying them by $\alpha_i$ (see (\ref{matv1}) and (\ref{matv2}) for details).

One can rearrange the columns of matrix $B$ in such a way that the updated matrix contains zero-blocks in the top-right corner (matrix $M$). This can be achieved by splitting the columns of $B_1$ and $B_2$ into sub-block, where the first sub-block contains $m+1$ columns and all the following sub-blocks contain $m$ columns (sub-blocks $B^{1,1},\dots,B^{1,n}$ and $B^{2,1},\dots,B^{2,n}$). Then 
$$M=[B^{11},B^{21},B^{21},B^{22},\dots,B^{1n},B^{2n}].$$
Note that $\rank(B)=\rank(M).$ 

According to the numerical experiments~\cite{Four_Models}, $B$ is a rank-deficient matrix and therefore $B^TB$ is  singular.  As a consequence, the normal equations  method is not efficient and therefore, more robust methods such as QR decomposition or SVD are required to solve this optimisation problem~\cite{SL, ABLSP, TB}. Note that these methods are substantially more expensive than normal equations.

\section{Singularity study}\label{sec:singularity_study}
Let us point out that the matrices $B,~B_1$ and $B_2$ can be expresses as a block lower triangular matrix $M$ such that
 \be\label{mat2}
M= \left[
 \begin{array}{cccccccccc}
 A_{11} & 0 &  \cdots & 0 & 0 \\
 A_{21} & A_{22} & 0 & \cdots& 0  \\
 \vdots & \vdots& A_{33} & 0 & 0  \\
 \vdots & \vdots & \vdots& \ddots & \vdots\\
 A_{n1} & A_{n2}& A_{n3}& \cdots & A_{nn}
 \end{array}
 \right ]
 \,.\ee
 where $A_{j1}$, $j=1, \dots,n$ has $N/n$ rows and $m+1$ columns, $A_{jk}$, $j,k=2, \dots,n$, $j\ge k$ has $N/n$ rows and $m$ columns and the top-right corner ($A_{jk}\,,$ for $k>j$, $j,k=1, \dots, n$) contains zeros since $$\max\{0,t_i-\theta_{l-1}\}=0,~\text {for~all}~t_i\leq\theta_{l-1}, ~l=2,\dots,n.$$

In some cases, one or more of the time moments can coincide with the corresponding spline knots.  To avoid possible ambiguity, we assign the time moments into the subintervals according to the following subdivision:
  \be\label{eq:int_division}[\theta_0,\theta_1],~(\theta_{k-1},\theta_k],~k=2,\dots n.
  \ee
Therefore, the first subinterval includes both borders while the other subintervals only include the right border.

Note that all the diagonal blocks of~(\ref{mat2})  are rectangular matrices ($A_{ii},~i=1,\dots,n$). The number of columns in $A_{11}$ is $m+1,$ while the number of columns in $A_{ii},~i=2,\dots,n$ is $m$. The number of rows in $A_{ii},~i=1,\dots,n$ coincided with the number of time moments assigned to the $i-$th subinterval.

\subsection{Model 1: singularity study}\label{subsec:singularity_study_model1}

In this section, we develop a sufficient condition for non-singularity of Model~1 (oscillation around ``zero''). If this condition is satisfied, we can guarantee that the corresponding matrices are non-singular and therefore, one can apply the normal equations method that is well-known to be fast and efficient for such problems. The following theorem holds. 

\begin{thm}\label{thm:model1}(sufficient non-singularty condition). Suppose that the spline degree is $m$, the number of subintervals is $n$ and the corresponding spline knots are 
$$\theta_0=t_1\leq\theta_1\leq\dots\leq\theta_{n-1}\leq\theta_n=t_N\,.$$
Matrix $B$ is non-singular when the following inequalities satisfy
\begin{equation}\label{eq:thm_model1}
N_1-Z_1\geq m+1~\text{and }~ N_k-Z_k\geq m,~k=2,\dots, n,
\end{equation}
where $N_{k},~k=1,\dots, n$ is the total number of recordings for $k-$th subinterval established in~(\ref{eq:int_division}) and $Z_k,$ $k=1,\dots,n$ is the number of time moments $t_i$ in the $k-$th subinterval, such that $g(t_i)=0$.
\end{thm}
{\bf Proof:}   Our proof is based on  two important facts.
\begin{enumerate}
\item For a rectangular matrix,  changing the order of the rows or multiplying a row by non-zero constants (elementary row operations) do not change the rank of the matrix.
\item  The determinant of a square Vandermonde matrix
 \begin{align}\label{van1}
 V=\left(
 \begin{array}{cccc}
 1 & X_{1} & X_{1}^2 &\dots X_{1}^{m} \\
 1 & X_{2} & X_{2}^2 &\dots X_{2}^{m} \\
 1 & X_{3} & X_{3}^2 &\dots X_{3}^{m}\\
 \vdots&\vdots&\vdots&\vdots\\
 1 & X_{m+1} & X_{m+1}^2 &\dots X_{m+1}^{m}
 \end{array}
 \right )\,,
 \end{align}
 can be expressed as 
 \begin{align*}
 \det(V)=\prod_{1\le i< j\le m}(X_{j}-X_{i})\,,
 \end{align*}
 and therefore, it  can not be zero if all $X_i,~i=1,\dots,m$ are distinct.
\end{enumerate}
Note that the rows of the diagonal blocks $A_{11}$  of~(\ref{mat2}) are Vandermonde matrix rows that are multiplied by a constant 
$\alpha_j,~j=1,\dots,N_1,$ where $N_1$ is the number of time moments assigned to the first interval. Therefore, if $N_1-Z_1\geq m+1$ then, it is possible to extract $(m+1)$ linearly independent rows from the first $N_1$ rows of $B$. 

For the second interval the situation is similar, but each row is multiplied by 
$$\alpha_j\times (t_j-\theta_1),~j=N_1+1,\dots, N_1+N_2.$$ 
Since for the second interval none of the time moments can coincide with $\theta_1$  one can conclude that 
$$(t_j-\theta_1)>0,~j=N_1+1,\dots, N_1+N_2.$$ 
Therefore,  if $N_2-Z_2\geq m$ then it is possible to extract $m$ linearly independent rows from the block of the rows $N_1+1,\dots, N_1+N_2$ of $B$. These rows will be also linearly independent with the $m+1$ rows extracted from the first $N_1$ rows of $B$.  

By continuing the process, we finally have $mn+1$ linearly independent rows and therefore, matrix $B$ is indeed full-rank. \qed


Generally, it is not always easy to estimate $Z_k,~k=1,\dots, n$. In Section~\ref{sec:applications} we give an example where $g(t)$ is a periodical (sine) function and therefore, the corresponding $Z_k,~k=1,\dots,n$ can be estimated through the corresponding frequencies.


\subsection{Model 2}\label{subsec:singularity_study_model2}



Consider $A_{11}$ and divide it into two parts: any $m+1$ rows form the bottom part and the remaining rows form the top part. Therefore, $A_{11}$ contains four sub-blocks $B_{111},~B_{121},~B_{211}$  and $B_{221}.$ The last index indicates that we are working with the block $A_{11}$ (first interval). Therefore,
\begin{equation}
A_{11}=\left(
\begin{array}{cc}
B_{111}&B_{121}\\
B_{211}&B_{221}\\
\end{array}
\right).
\end{equation}

 $B_{221}$ is a full-rank block in $A_{11}$ (Vandermonde-type matrix) and if we apply equivalent row operations one can obtain zeros everywhere in sub-block $B_{121}$ using the last $m+1$ rows of $A_{11}.$  Therefore, there exists a unique set of $\lambda^1_{ij},~j=1,\dots, m+1,~i=1,\dots, m+1,$ such that

\begin{equation}\label{rowsB121}
B_{121}^j=\sum_{i=1}^{m+1}\lambda^1_{ij}B_{221}^i\,,
\end{equation}
where $B_{121}^j$ is the $j-$th row of $B_{121}$ and $B_{221}^i$ is the $i-$th row of $B_{221}$.

The rows of $B_{111}$ and $B_{211}$ denote by $B_{111}^j$ and $B_{211}^j$ for ~$j=1,\dots,m+1$ respectively. Then, the rows $\tilde{B}_{111}^i,~i=1,\dots,m+1$ of the updated block $\tilde{B}_{111}$ (obtained from $B_{111}$ by equivalent row operations) are as follows:
\begin{equation}\label{eq:rows}
\tilde{B}_{111}^j=B_{111}^j-\sum_{i=1}^{m+1}\lambda^1_{ij}B_{211}^i.
\end{equation}
Therefore,
\begin{equation}\label{eq:row_to_matrix}
\tilde{B}_{111}=B_{111}-\Lambda_1^T B_{211},
\end{equation}
where $\Lambda_{1}$ is an $(m+1)\times(m+1)$ matrix with the $\{ij\}$-th element equals $\lambda^1_{ij}.$
Note that $A_{11}$ is full-rank if and only if
\begin{equation}\label{eq:tildeB}
\tilde{B}_{111}=B_{111}-\Lambda^T_{1} B_{211}\,,
\end{equation}  
is full-rank. Similar reasoning is applied to remaining intervals. Taking into account that 
\begin{itemize}
\item the dimension of each block $A_{ii},~i=2,\dots,n$ is $N_i\times m\,,$ where $N_i$ is the number of time moments assigned to the $i-$th interval and
\item the left border point is not included in any of each interval.
\end{itemize}

In general, for $1<i\leq n$, (\ref{eq:tildeB}) can be rewritten as follows
\begin{equation}\label{eq:tildeBg}
\tilde{B}_{11i}=B_{11i}-\Lambda^T_{i} B_{21i},
\end{equation}
 where  $\Lambda_{i}$ is an $m\times m$ matrix with the $\{jk\}$-th element equals $\lambda^i_{jk}.$ Hence, the following theorem holds (sufficient condition for non-singularity of Model~2).
%

\begin{thm}\label{thm:model2}
Suppose that the spline degree is $m$, the number of subintervals is $n$ and the corresponding spline knots are 
$$\theta_0=t_1\leq\theta_1\leq\dots\leq\theta_{n-1}\leq\theta_n=t_N\,.$$
It is possible to construct the bottom blocks in each matrix $A_{ii}$ such that the corresponding matrices in (\ref{eq:tildeBg}) are full-rank then, $B$ is full-rank too.
\end{thm}

Note that if $g(t)$ is constant then, the condition of Theorem~\ref{thm:model2} are not satisfied. Moreover, it indicates that the corresponding matrix in Model 2 is rank-deficit since there exists a column that is obtained by multiplication of another column by $a=g(t_i)$, $i=1,\dots,n.$



\section{Application to signal processing}\label{sec:applications}
In this section, we give an example of how our conditions can be applied in a particular problem in signal processing. In addition, we propose an algorithm for signal approximation where the choice of optimisation techniques is based on the corresponding singularity study.  
\subsection{Model 1}\label{subsec:application_model_1}

An EEG (electroencephalogram, also known as brain wave) signal is modeled as a sine wave 
\be\label{wave1} W_{1}=S_{m}(\vec x, \boldsymbol{\tha}, t)  \sin(\om t + \shi)\,,\ee
where $S_{m}$ is the spline function defined in (\ref{eqn4}) whose $\boldsymbol{\tha}=(\theta_{1},\dots,\theta_{n-1})$ are equidistant therefore, for each combination of $\om$  and $\shi$ the corresponding optimisation problem is
\be\label{model1}
\min_{\vec x} \sum_{i=1}^N (y_{i}-S_{m}(\vec x, \boldsymbol{\tha}, t_{i}) \sin(\om t_{i}+\shi))^2.
\ee

This is an LLSP. To achieve the best combination of  $\om$  and $\shi$ we run a double loop over the defined intervals for $\omega$ and $\tau$,  keeping the optimal combinations, that is combinations with the lowest objective function values (see section~\ref{subsec:algorithm} for details). Optimisation problem (\ref{model1}) is reformulated as 
\be
\min_{\vec x} \sum_{i=1}^N (y_{i}-M\vec x^T)^2\,, \quad\mbox{or} \quad\min_{\vec x} \|M\vec x-y_{i}\|_{2}^2\,,
\ee
where $y_{i}\in \R^N$ are the recorded signals at $t_{i}\in \R^N$, $\vec x\in \R^{mn+1}$ and $M$ is a matrix with $N$ rows and $mn+1$ columns of the form $S_{m}(\vec x, \boldsymbol{\tha}, t_{i}) \sin(\om t_{i}+\shi)$. If $M\in \R^{N\times mn+1}$ is a full-rank matrix then, this LLSP can be solved through the systems of normal equations. 

In this application, there are natural restrictions on $\omega$ and $n$ that are considered in~\cite{ANZ_us}, 
\begin{itemize}
\item Frequency $\omega$ is a parameter that is normally assigned by a manual scorer (medical doctor). Therefore, the value for frequency is bounded from above (by 16~Hz) and restricted to integer (due to human scorer's perception limitations).  
\item The duration of the events is between 0.5 and 3~seconds (these events are called K-complexes)~\cite{ANZ_us}. It is not reasonable to consider any  interval shorter than 1~second. Since the duration of the original signal is 10~seconds therefore, the number of subintervals $n$ can not exceed 10.
\end{itemize}
   
   We need to show that there are $mn+1$ linearly independent rows and therefore, $B$ is a full-rank matrix (see Theorem~\ref{thm:model1}). There are three main reasons why there may be fewer than $mn+1$ linearly independent rows.
   \begin{enumerate}
   \item Some of $\sin(\omega t_i+\tau)$ are zero then, in this case, zero rows will appear. This can happen when the frequency is high. 
   \item There are too many intervals $n$, since in this case $mn+1$ is large. 
   \item The degree $m$ of the corresponding polynomial is high so, in this case, $mn+1$ is large. There is no (application based) upper bound for $m$ therefore, the choice of $m$ is based on the complexity of the signal and computer resources.
   \end{enumerate}
 
Before we start proving that $M$ is a full-rank matrix, we have to estimate the number of constants $\sin(\omega t_i+\tau)$ being zero (estimation of $Z_k,~k=1,\dots,n$).   It can happen no more than $2\om D +1$ times for a  $D$ seconds duration of an EEG. Suppose that the knots $\theta_{1},\dots,\theta_{n}$ (switches from one polynomial to another) are equidistant. According to the reported experiments in~\cite{ANZ_us} $N=1000$, $n=5$, $m=4$ and $\omega$ did not exceed 16~Hz. The duration of each signal segment is 10~seconds therefore, the duration of each subinterval is $$D/n=2~\text{seconds}\,.$$ Then,
$$N_k=200,~Z_k\leq 2\omega\times D/n+1=65\,,~k=1,\dots,n\,,$$
and therefore,
$$N_k-Z_k=200-65\geq m=4,~k=1,\dots, n.$$
Hence, due to the Theorem~\ref{thm:model1} matrix $M$ is non-singular and therefore, the normal equations method can be applied. 

\subsection{Model 2}\label{subsec:application_model_2}
Similar to Model~1, the knots are equidistant, $N=1000$, $n=5$, $m=4$ and $\omega$ did not exceed 16~Hz. An EEG  signal is modeled as a sine wave that is shifted vertically by a spline function as
\be\label{shiftwave1}
W_{2}=S_{m}(\vec x_{1}\,, \boldsymbol{\tha}\,, t) \sin(\om t+\shi)+S_{m}(\vec x_{2}\,, \boldsymbol{\tha}\,, t)\,,
\ee
therefore, the corresponding optimisation problem is 
\be\label{model2}
\min_{\vec x} \sum_{i=1}^N (y_{i}-S_{m}(\vec x_{1}\,, \boldsymbol{\tha}\,, t_{i}) \sin(\om t_{i}+\shi)-S_{m}(\vec x_{2}\,, \boldsymbol{\tha}\,, t_{i}))^2\,,
\ee
where $y_{i}\in \R^N$ are the recorded signal at $t_{i}$ for $i=1,2,\dots,N$ and $\vec x=[\vec x_{1}; \vec  x_{2}]$. The dimension of this problem is $2mn+2$. The optimization problem (\ref{model2}) can be rewritten as
\be
\min_{\vec x} \sum_{i=1}^N (y_{i}-B\vec x^T)^2\,, \quad\mbox{or} \quad\min_{\vec x} ||B\vec x-y_{i}||_{2}^2\,,
\ee
where $B\in \R^{N\times (2mn+2)}$.  Matrix $B$ in the optimisation problem (\ref{model2}) is 
$$B^{N \times (2mn+2)}=[B_{1}^{N \times mn+1} \quad  B_{2}^{N \times mn+1}],$$ where $B_{1}^{N\times mn+1}$ and $B_{2}^{N \times mn+1}$ are detailed in~(\ref{matv1}) and (\ref{matv2}). 

Numerical experiments with Model~2 (see~\cite{Four_Models}) indicate that $B$ is a rank-deficient matrix. This can be anticipated, since the conditions of Theorem~\ref{thm:model2} are not satisfied. In this case, instead of checking the rank of an $N\times (2mn+2)$ matrix one needs to check the rank of several  $N/n\times (mn+1)$ matrices. Since the number of such matrices can be large, in most cases it is more efficient to check the singularity of the original matrix. Numerical experiments in~\cite{Four_Models} show that the corresponding matrices are not full-rank and therefore, an SVD-based method is applied instead of normal equations. 

 \subsection{Algorithm implementation}\label{subsec:algorithm}
 In this section, we present an algorithm for solving~(\ref{model1}) and~(\ref{model2}). In most practical problems, $\omega$ and $\tau$ are not known in advance and therefore, there should be a procedure for choosing them. One way  is to consider them as additional variables and optimise them. This approach is not very efficient since the corresponding optimisation problems become non-convex and can not be solved fast and accurately~\cite{ANZ}. Therefore, we can assign exact values from defined intervals of $\om$ and $\shi$ that form a fine grid (using double loops) instead of optimising them directly. Then, we solve the corresponding LLSPs  and keep the best obtained results~\cite{ANZ_us}. 
  
 The following algorithm can be used to solve a sequence of LLSPs. In this algorithm, $\omega_0$ and $\omega_f$ are the initial and final value for $\omega$. Similarly, $\tau_0$ and $\tau_f$ are the initial and final values for $\tau$.
\begin{center}
 {\bf Algorithm 1: Signal approximation through LLSPs}
 \end{center}

 \begin{algorithmic}[1]
 \State {Specify the initial and final values for the frequency ($\omega_0$ and $\omega_f$) and shift ($\tau_0$ and $\tau_f$)}
    
 \For {$\om=\om_{0}:\om_{f}$}
 \For {$\shi=\shi_{0}:\shi_{f}$}

 \State{Solve the corresponding optimisation problem (LLSP) with fixed $\omega$ and $\tau$;  and record the minimal value of the objective function.}
 \EndFor
 \EndFor
 \end{algorithmic}

In this algorithm, the normal equations method can be used for solving LLSPs with non-singular matrix while QR decomposition or SVD should be used for singular cases.

\section{Conclusions and further research directions}\label{sec:conclusions}
Most linear least square problems can be solved using the system of normal equations if the corresponding matrix is non-singular otherwise one needs to apply a more robust (and time-consuming) approach (e.g., QR decomposition and SVD). We consider two types of linear least squares problems and develop a procedure which enables us to identify when the corresponding matrix is non-singular. Basing on the outcomes of this procedure, one can choose a more efficient method for solving the corresponding linear least squares problems.

Currently, we consider three main future research directions.
\begin{enumerate}\item The development of necessary and sufficient conditions for non-singularity verification. 
\item The development of more flexible models where the spline of vertical shift does not have the same degree and knots location as the main spline  (multiplied by  prototype functions). 
\item
The extension of the results to the case when other types of functions (not necessary polynomial splines) are used to construct the corresponding approximations. 
\end{enumerate}

\bibliographystyle{plain}
\bibliography{mybib}

\end{document}